\documentclass[12pt]{amsart}
\usepackage{amscd}
\usepackage{verbatim}

\usepackage[dvips]{graphics}
\usepackage[cp866]{inputenc}
\usepackage{graphicx}
\usepackage{epsfig}

\usepackage{amsmath,amssymb,amscd,}
\usepackage[all]{xy}

\usepackage{amssymb}

\newtheorem{thm}{Theorem}[section]

\newtheorem{prop}[thm]{Proposition}

\newtheorem{ex}[thm]{Example}

\newtheorem{lem}[thm]{Lemma}

\newtheorem{conq}[thm]{Question}
\theoremstyle{definition}

\newtheorem{defn}[thm]{Definition}

\title[Pseudo-tilings and hyperbolic virtual polytopes]
{Planar pseudo-triangulations, spherical pseudo-tilings  and
hyperbolic virtual polytopes}

\author{Gaiane Panina}


 \email{panina@iias.spb.su}

\begin{document}
\begin{abstract}

We wish to draw attention to an interesting and promising
interaction of two theories.

On the one hand, it is the theory of \textbf{pseudo-triangulations}
which was useful for implicit solution of the
 carpenter's rule problem  and proved later to give a nice tool for
 graph embeddings.

On the other hand, it is the theory of \textbf{hyperbolic virtual
polytopes} which arose from an old uniqueness conjecture for convex
bodies (A. D. Alexandrov's problem): suppose that  a constant $C$
separates (non-strictly) everywhere  the principal curvature radii
of a smooth 3-dimensional convex body $K$. Then  $K$ is necessarily
a ball of radius $C$.

The two key ideas are:
\begin{itemize}
    \item Passing from planar pseudo-triangulations to spherical pseudo-tilings,
    we avoid non-poited vertices.
Instead, we use pseudo-di-gons. A theorem on spherically embedded
Laman-plus-one graphs is announced.
    \item The difficult problem of
 hyperbolic polytopes constructing  can  be reduced to
finding spherically  embedded graphs.
\end{itemize}

\end{abstract}

\maketitle

\setcounter{section}{-1}

\section{Introduction}

Pseudo-triangulations are opposite to the traditional planar graph
embeddings - they are as non-convex as possible.

As a parallel phenomenon, hyperbolic virtual polytopes are opposite
to convex polytopes: convexity is replaced by saddle property.

The two theories have a nice interaction, which is demonstrated in
the paper.

Even at first glance  one can see  that
 pseudo-triangulations look very much like the fans of
 hyperbolic virtual polytopes. Indeed, in both cases we have a pointed tiling.
 But whereas pseudo-triangulations  are planar drawings,
 the fans of
 hyperbolic virtual polytopes are spherical ones.

 However, the relationship is much deeper. As is shown in the paper,
 some items are absolutely
the same, some are easily adjustable, but some are different.

In the first two sections, we sketch very briefly the two theories,
omitting all details and applications and referring the reader to
detailed papers from the list of references.

In the third section, we bring the two theories together.

We show that for spherical embeddings, the usage of pseudo di-gons
allows pointed embeddings not only for Laman graphs, but also of
graphs with a greater number of  edges.
 Even some Laman-plus-$k$ graphs
(for any natural number $k$) admit a pointed embedding (see Example
3.7 for $k=5$).
 In particular, we announce  the following  theorem
on Laman-plus-one graph embedding. It is parallel to the results of
[7],[25] (quotated  also in section 1).

\bigskip

\textbf{Theorem on spherical embedding of Laman-plus-one graphs.}
(Theorem 3.3)

\begin{itemize}
    \item Each Laman-plus-one graph admits a straightened embedding
(all the edges are geodesic segments) in the 2-dimensional sphere
$S^2$ such that
\begin{enumerate}
    \item it generates a pointed pseudo-tiling of the sphere.
    \item The tiles are either pseudo-triangles or pseudo-di-gons.
    \item  The number of pseudo-di-gons equals 4.
    \item The embedding can be constructed inductively, via
    geometric Henneberg constructions starting from the fan of a hyperbolic tetrahedron
    (i.e. from the pointer embedding of $K_4$, see Fig.  1)
\end{enumerate}
    \item Any spherical embedding of a Laman-plus-one graph as a nice pseudo-tiling
    has exactly 4 pseudo di-gons.
    \item If a graph admits a straightened embedding   possessing the above properties 1-3,
    then it is a Laman-plus-one graph.\qed
\end{itemize}

Some  very natural questions on spherical pointed tilings are
formulated.

A simple example (Example  3.9) demonstrates how the methods of
pseudo-triangulations work for the sake of the hyperbolic polytopes,
and visa versa.

At the end of Section 4, we arrange the parallel terms from both the
theories in a kind of a dictionary. The correspondence of the
objects is not always straightforward,
 but we find better to skip additional technicalities
in order to stress similarity of the ideas.

\textbf{Acknowledgments }

The research was inspired by talks and discussions at the Ervin
Scr\"{o}dinger Institute Program "Rigidity and Flexibility" (2006)
organized by V. Alexandrov, I. Sabitov, and H. Stachel.

The author is also greatful to Nikolai Mnew for stimulating
conversations.

\section{Pseudo-triangulations}

\textbf{Topological preliminaries.}

In the paper, we consider only planar graphs.

Let $\Gamma$ be a graph with $v$ vertices and $e$ edges.

$\Gamma$ is a \textit{Laman graph} if $e=2v-3$ and every subset of
$k$ vertices spans at most $2k-3$ edges.

Laman graphs are interesting because they are\textit{ minimally
rigid graphs} - their generic embeddings (as     bar-and-joint
frameworks) are rigid.

$\Gamma$ is \textit{Laman-plus-one graph} (respectively,
\textit{Laman-plus-$k$}) if there is an edge (respectively, $k$
edges) such that after its removing the graph becomes a Laman graph.

Minimal  Laman-plus-one graphs are called  \textit{rigidity
circuits}. A generic embedding of a rigidity circuit possesses a
unique (up to a constant)  non-vanishing
 self-stress. Therefore, the embedded graph has a 3D lift (i.e. it can be represented
as a projection of a spatial polytope). Besides, such a graph has a
unique (up to a homothety)
 \textit{Maxwell reciprocal} (a 3-dimensional  dual to the graph polytope
 constructed by the stress).

Laman and Laman-plus-one graphs admit an inductive construction
starting with elementary graphs.
 At each step, a new vertex is added via one of the
following two \textit{Henneberg constructions}.

\begin{itemize}
    \item Henneberg 1 construction: add a new vertex,
    connecting it via two new edges to two old vertices.
    \item Henneberg 2 construction:  add a new vertex inside an
    old edge (and thus split the edge into two new ones),
     connecting the new vertex  by another new edge  to a third  vertex.
\end{itemize}

\begin{lem} (see [7], [28])
\begin{enumerate}
    \item A graph is Laman if and only if it admits
     an inductive construction starting with a graph with two vertices
and one edge.
    \item A graph is Laman-plus-one if and only if it admits
     an inductive construction starting with $K_4$ (the complete graph with 4 vertices).

     In both cases, each step is a Henneberg construction.\qed
\end{enumerate}
\end{lem}

\bigskip

\textbf{Geometric realizations.}

Consider a planar embedding of a graph $\Gamma$. We say that its
vertex is \textit{pointed} if one of the adjacent angles is greater
than $\pi$.

An  embedding is \textit{pointed} if all its vertices are pointed.

A \textit{pseudo-triangle} is a simple polygon (a non-crossing
planar broken line) which has exactly 3 convex vertices.

A \textit{pseudo-triangulation} is a tiling of a convex polygon in
the plane such that each tile is a pseudo-triangle.

A pseudo-triangulation is \textit{pointed} if all its vertices are
pointed.

A pseudo-triangulation is \textit{pointed-plus-one} if all its
vertices, except for exactly one vertex,
 are pointed.

\begin{thm} [7] (On planar embedding of Laman and Laman-plus-one  graphs)

\begin{itemize}
    \item A graph $\Gamma$ is a planar Laman graph if and only if
    it can be embedded in the plane as a pointed pseudo-triangulation.
    \item A graph $\Gamma$ is a planar Laman-plus-one  graph if and only if
    it can be embedded in the plane as a pointed-plus-one  pseudo-triangulation.

    In both cases the embeddings can be constructed inductively
    such that

    \end{itemize}
    \begin{itemize}
        \item Construction starts by an embedding of one-edge-graph (for Laman graphs)
        or by an embedding of $K_4$ (for Laman-plus-one graphs)
        \item On each step, we get a  pointed pseudo-triangulation
        (respectively, a pointed-plus-one  pseudo-triangulation)
        \item Each step is a geometric realization of a Henneberg construction.
        The construction is local: it does not change the positions of old vertices.\qed

    \end{itemize}

\end{thm}

\section{Hyperbolic virtual polytopes}

Roughly speaking, virtual polytopes are geometric realizations of
Minkowski difference of convex polytopes.

They were introduced originally by A. Pukhlikov and A. Khovanskij in
[10], appeared also in a different disguise in the polytope algebra
of P. McMullen [14].

Hyperbolic virtual polytopes [18-20] are virtual polytopes with
special saddle properties.

In the section, we try to give  a shortcut to necessary  notions.

\bigskip

Convex polytopes in $ \mathbb{R}^3$ form a semigroup $\mathcal{P}$
with respect to the Minkowski addition $\otimes$.

The semigroup  $\mathcal{P}$  is isomorphic to the semigroup of
 continuous convex piecewise linear (with respect to a \textit{fan}) functions
 defined on $ \mathbb{R}^3$.

 The isomorphism maps a convex polytope to its support function.

 \bigskip

 (A necessary reminding:  the support function of a polytope
 is piecewise linear with respect to
 some conical tiling of the space. To visualize the tiling, we intersect it with a unite
 sphere centered at $O$ and get a spherical \textit{fan} of the polytope.
 It is a spherical tilings, all tiles are convex.
 In some sense, a polytope can be considered as the Maxwell's reciprocal of its fan.)

 \bigskip

 Passing to the Grothendieck group $\mathcal{P}^*$ (it is the group of formal differences
 of convex polytopes)
 which is called the \textit{group of virtual polytopes}, only
 the convexity property disappears.
 Thus we get a group isomorphism

\bigskip
 \textbf{virtual polytope  $\longleftrightarrow$  continuous  piecewise linear
 (with respect to a
 fan) function
 defined on $ \mathbb{R}^3$.}

\bigskip

The skeleton of the fan has a self-stress.
 A virtual polytope (it can be considered as  the Maxwell's reciprocal of its fan)
 can be represented geometrically as a polytopal function [14]
 or as a closed polytopal surfaces [20].

We don't mind (and can not avoid) self-crossing 3D reciprocals.
 This is because
hyperbolic polytopes (considered as spatial piecewise linear
surfaces) are necessarily self-crossing (except for degenerated
cases as hyperblic tetrahedron), see Example  2.3.

 Given a virtual polytopes, the tiles of its fan can be non-convex.

Recall that the support function of a convex polytope is convex,
i.e. its graph is a convex surface  (it is reasonable to consider
either the spherical graph or
 the collection of affine graphs  [18],[19]).

Among virtual polytopes  we point out the class  of hyperbolic
virtual polytopes.

\begin{defn}
A virtual polytope is hyperbolic if the graph of its support
function is a saddle surface.
\end{defn}

This definition arose quite natural from the following conjecture.

\textit{Given a smooth 3-dimensional convex body $K$ and a constant
$C$ such that    $R_1 \leq C \leq R_2$ holds at each point of
$\partial K$ ($R_1$ and  $ R_2$ stand for the principal curvature
radii of $K$),
 the body $K$ is necessarily a ball of radius $C$.}

\bigskip

 The conjecture proved to be wrong (see [12], [18]), and  here is a way of constructing
 counterexamples (which are unexpectedly  diverse).
 \begin{itemize}
    \item Construct a hyperbolic polytope (this is the most difficult step, for hyperbolic
    polytopes are very rare phenomena among virtual polytopes)
    \item Smoothen its support function $h$ (preserving saddle property)
    \item Add to $h$ the support function of a ball (which is sufficiently large to
    make
    the sum convex). The result is the support function of a counterexample
    to the conjecture.
 \end{itemize}

In the framework of the theory of hyperbolic polytopes, pointed
spherical tilings  appear due to the following simple observation.

\begin{lem} (see [20])
 \begin{itemize}
    \item The fan of a virtual polytope $K$ is a  pointed  tiling $\Rightarrow$ $K$ is hyperbolic.
    \item If $K$ is simplicial, then

    the fan of $K$ is a  pointed  tiling $\Leftrightarrow$ $K$ is hyperbolic.\qed
 \end{itemize}
\end{lem}

\begin{ex}

Figure 1 presents  the fan of  the \textit{hyperbolic tetrahedron}.
It is the simplest hyperbolic polytope.

The hyperbolic tetrahedron is useless for the above conjecture (for
this polytope, the smoothing technique of [18] does not work), but
it proved to be the starting point for Laman-plus-one graphs
embeddings (Theorem 3.3).
\end{ex}

\pagebreak

\begin{figure}

\centering
\includegraphics[width=10 cm]{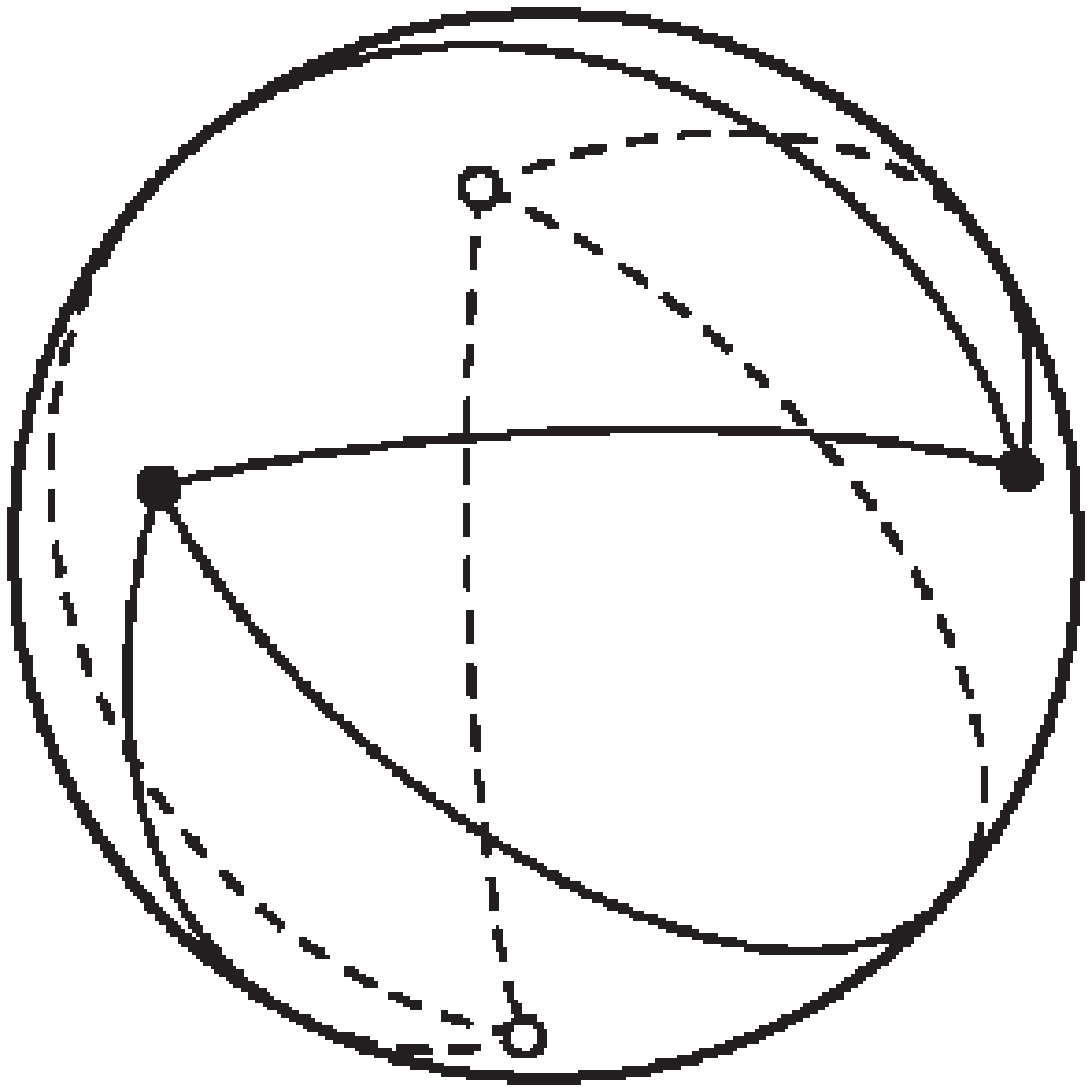}
\caption{}\label{}

\end{figure}

\begin{ex}

Figure  2 depicts a hyperbolic  polytope (viewed as a complicated
self-intersecting 3D surface)  and its fan.

The hyperbolic polytope has 8 \textit{horns} - the non-saddle
vertices. By duality, they correspond to 8 pseudo di-gons (see the
definition in Section 3).

The di-gons are marked grey. Note that only half of each of di-gons
is visible.

For each hyperbolic polytope,   \textbf{horns are   dual to
pseudo-di-gons}. For a simplicial hyperbolic polytope, duality maps
bijectively horns of the polytope to the pseudo-di-gons of its fan
[20].
\end{ex}

\pagebreak

 \begin{figure}

\centering
\includegraphics[width=14 cm]{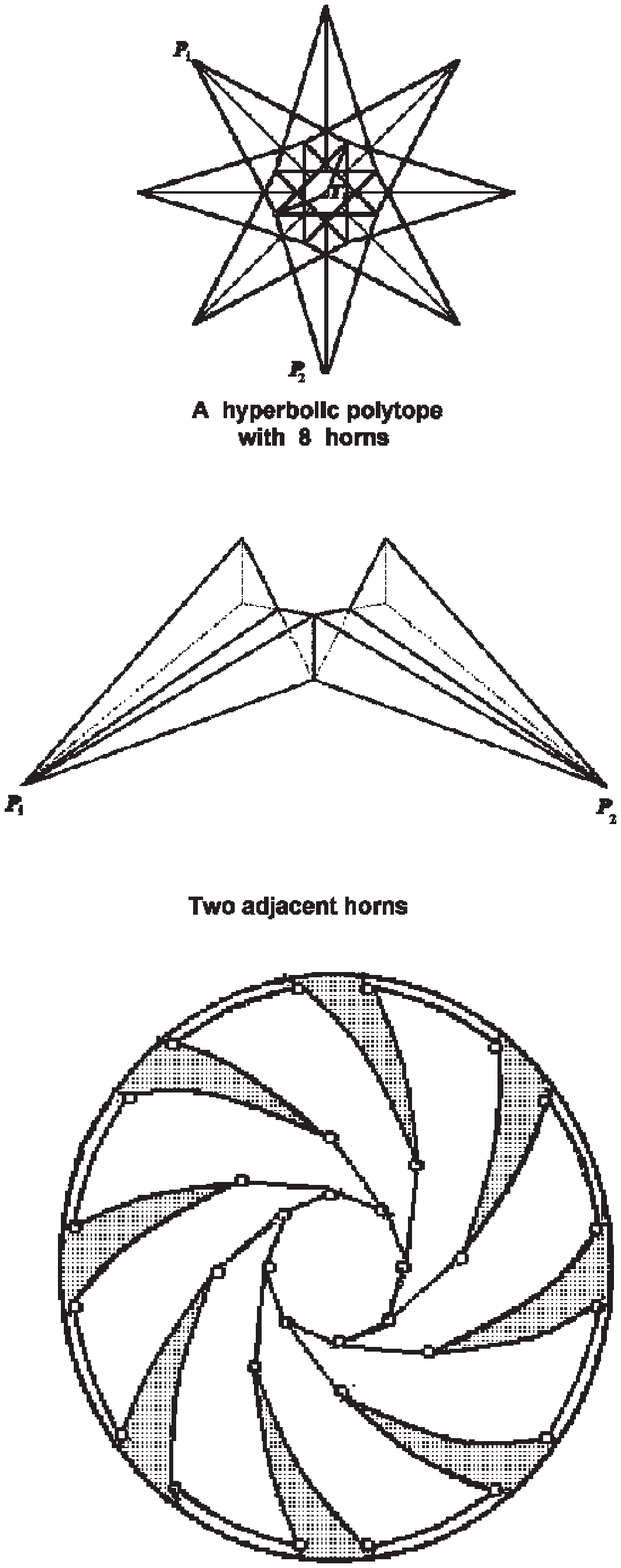}
\caption{}\label{}

\end{figure}

It makes sense to color the edges of the fan of a hyperbolic
polytope  $K$ red and blue. The support function of $K$ is concave
up along the red edges and concave down on the blue ones.

 The theory of hyperbolic polytopes has another nice applications. Here we list
 some problems from classical geometry which have been solved using this theory.

 \begin{enumerate}
    \item A refinement of A.D. Alexadrov's unieqness theorem for 3D-polytopes
    with non-insertable pairs of parallel faces [19].
    \item Extrinsic geometry of saddle surfaces with injective Gaussian mapping [18-21].
    \item Isotopy problem for saddle surfaces  [21].
 \end{enumerate}

\section{Putting the pieces together}

>From now on, we consider graph embeddings in the sphere $S^2$.

\bigskip

\textbf{The first key idea is: for the sake of hyperbolicity, we
avoid non-poited vertices. Instead, we use pseudo-di-gons.}

\bigskip

A \textit{pseudo di-gon} is a spherical polygon  (a non-crossing
broken line embedded in the sphere with a fixed interior domain; all
segments are geodesic segments) which has exactly two convex
vertices.

\begin{defn}
A \textit{\textbf{nice pseudo-tiling }}is a spherical tiling which
is \begin{itemize}
    \item pointed
    \item each tile is either a (spherical) pseudo-triangle or a pseudo-di-gon.

\end{itemize}

\end{defn}

 The  di-gons of a nice pseudo-tiling are of a particular interest
 from the viewpoint of both theoties. Firstly, they correspond by duality to the
 horns of hyperbolic polytopes. Secondly, their number determins  the  Laman-type
 counts.

\begin{prop}
 The number of di-gons  can range from 0 to infinity.
\end{prop}

Proof. To construct
 an embedded graph with $k>3$ di-gons, take the fan of a hyperbolic polytope with $k$
 horns (see [20]) and find its nice pseudo subtiling.
 For $k<4$, it is easy.\qed


\begin{thm}
\begin{itemize}
    \item Each Laman-plus-one graph admits a straightened embedding
(all the edges are geodesic segments) in the 2-dimensional sphere
$S^2$ such that
\begin{enumerate}
    \item it generates a nice  pseudo-tiling of the sphere;
    \item there are exactly 4  pseudo di-gons;
    \item the embedding can be constructed inductively, via
    geometric Henneberg constructions starting from the fan of the hyperbolic tetrahedron
    (i.e. from a pointer embedding of $K_4$, see Fig.  1).
\end{enumerate}
   \item Any embedding of a Laman-plus-one graph as a nice pseudo-tiling
    has exactly 4 pseudo di-gons.

    \item If a graph admits a straightened embedding   possessing the above properties 1-2,
    then it is a Laman-plus-one graph.
\end{itemize}

\end{thm}


Proof  (a sketch).

\textbf{1.}   The proof repeats  that  of Theorem  3.2, [7].

Two items are essential.

On one hand, dislike  [7],  we do not care about geometric
realization of a combinatorial pseudo-triangulation (we do not
prescribe what angles should be greater than $\pi$).

On the other hand, the geometric shape of spherical pseudo-triangles
and pseudo-di-gons can be bad and can cause obstacles for applying
Henneberg constructions.

This motivates the following definition.

A pseudo-triangle (or a pseudo di-gon) is called \textit{H-good} if
it admits geometricaly any Henneberg construction.

The H-goodnes can be expressed in terms of \textit{feasible
regions}(see [7]) of the tile.

All planar pseudo-triangles (and therefore, all their spherical
images) are good (as is proven in [7]).

When applying a geometrical Henneberg construction, it is always
possible to preserve H-goodnes of all the tiles.

Thus we get an algorithm of the desired embedding (which is nearly
the same as in [7]):

\begin{itemize}
    \item By Lemma  1.1, we have an inductive topological  Henneberg
construction of the graph starting from $K_4$.
    \item Take the pointed embedding of $K_4$ (Fig.  1). It is H-nice.
    \item Realize geomerticaly step by step the Henneberg constructions,
    preserving H-goodness.
    Note that the number of pseudo di-gons does not change.
\end{itemize}

\textbf{2.} The proof repeats literally the corners counts from [7]
and recalls very much color changes counts for hyperbolic fans [20].

Denote by $v$ the number of vertices, by  $c$ the number of
\textit{corners } ( i.e. the number of angles which are greater than
$\pi$), by $e$  the number of edges, by $f_3$ the number of
pseudo-triangles, and by $f_2$ the number of pseudo di-gons.

We have $v-e+f_2+f_3 =2$ (Euler formula),

$e=2v-2$  (Laman-plus-one count),

$c=2f_2+3f_3$  (obvious),

and $ c=2e-v$ (corners count), which easily complete the proof.

\textbf{3.} Assume that an embedding of a graph $\Gamma$ generates a
nice pseudotiling with 4 pseudo-di-gons.

The above counts imply that $e=2v-2$.

Fix $k$ vertices and denote by  $\Gamma '$ the spanned subgraph. It
generates a pointed spherical tiling. Then the number of di-gons is
not greater than 4 because no di-gon
    admits a pointed subtiling into a collection of pseudo-triangles.
Similar counts complete the proof. \qed


\begin{conq}
What part of combinatorics of a Laman-plus-one graph embedding can
be prescribed (as is done in  [7])?
\end{conq}

\begin{conq}
Is there any canonical Laman graphs embedding in the sphere?

(Note that we have already at least two different  pointed spherical
embeddings for a Laman graph: the first one comes from its planar
pointed embedding raised to the sphere; to get the second one, just
add an edge, embed the result according to Theorem 3.3, and then
erase the edge.)
\end{conq}

\begin{conq}
There exist nice pseudo-tilings  with no di-gons (all tiles are
pseudo-triangles). What is a characterization of the set of planar
graphs admitting such an embedding?
\end{conq}

\begin{ex}
The tiling from Example 2.5 obviously admits a subtiling which is a
nice pseudo-tiling. It gives a pointed embedding of a graph with
$e=2v+2$.
\end{ex}

Thus manipulations with the number of di-gons enable us to embed
graphs with many edges.

\begin{conq}
     Is each  nice pseudotiling such that the number of pseudo-di-gons
    equals $3+k$   generated  by a Laman-plus-k graph?
\end{conq}

\bigskip

\textbf{Another key idea: the difficult problem of constructing of
hyperbolic polytopes (3D objects) can sometimes be reduced to
constructing embedded graphs  (2D objects).}

The following simple example demonstrates how it can work. Note that
the first item (which is already known, see [13] and [20])
 looks  quite trivial.
Anyhow, the statement was not so trivial and needed much efforts
three years ago.

\begin{ex}\begin{itemize}
    \item There exists a hyperbolic polytope with 4 horns.
    \item We present a Laman-plus-one-graph (a rigidity circuit)
 embedded in $S^2$ as a pointed pseudotiling with 4 pseudo-di-gons.

\end{itemize}
Proof.

The usual counts show that the  spherical graph in Figure 3  is a
rigidity circuit. It has a non-vanishing self-stress. The
self-stress gives a virtual polytope, which is hyperbolic because
the tiling is pointed.  The number of horns equals 4 because it
equals the number of di-gons (marked grey).\qed
\end{ex}

Obviously, the edge coloring of the fan of a hyperbolic virtual
polytopes (= of an embedded self-stressed   graph) reflects the sign
pattern of the stress.

\begin{conq}
     Given an embedded rigidity circuit, is it possible to detect
     the sign pattern of its (unique) self-stress by the combinatorics of the embedding
     (i.e. by corners information)?
\end{conq}

\begin{conq}
     There seems to be no straightforward spherical analog for \textit{expansive motions}
     [25]. Does there exist a parallel statement for spherically embedded graphs
     which exploits  the same underlying reasons (duality  together with
     mountain-valley arguments)?
\end{conq}

\pagebreak
\begin{figure}

\centering
\includegraphics[width=10 cm]{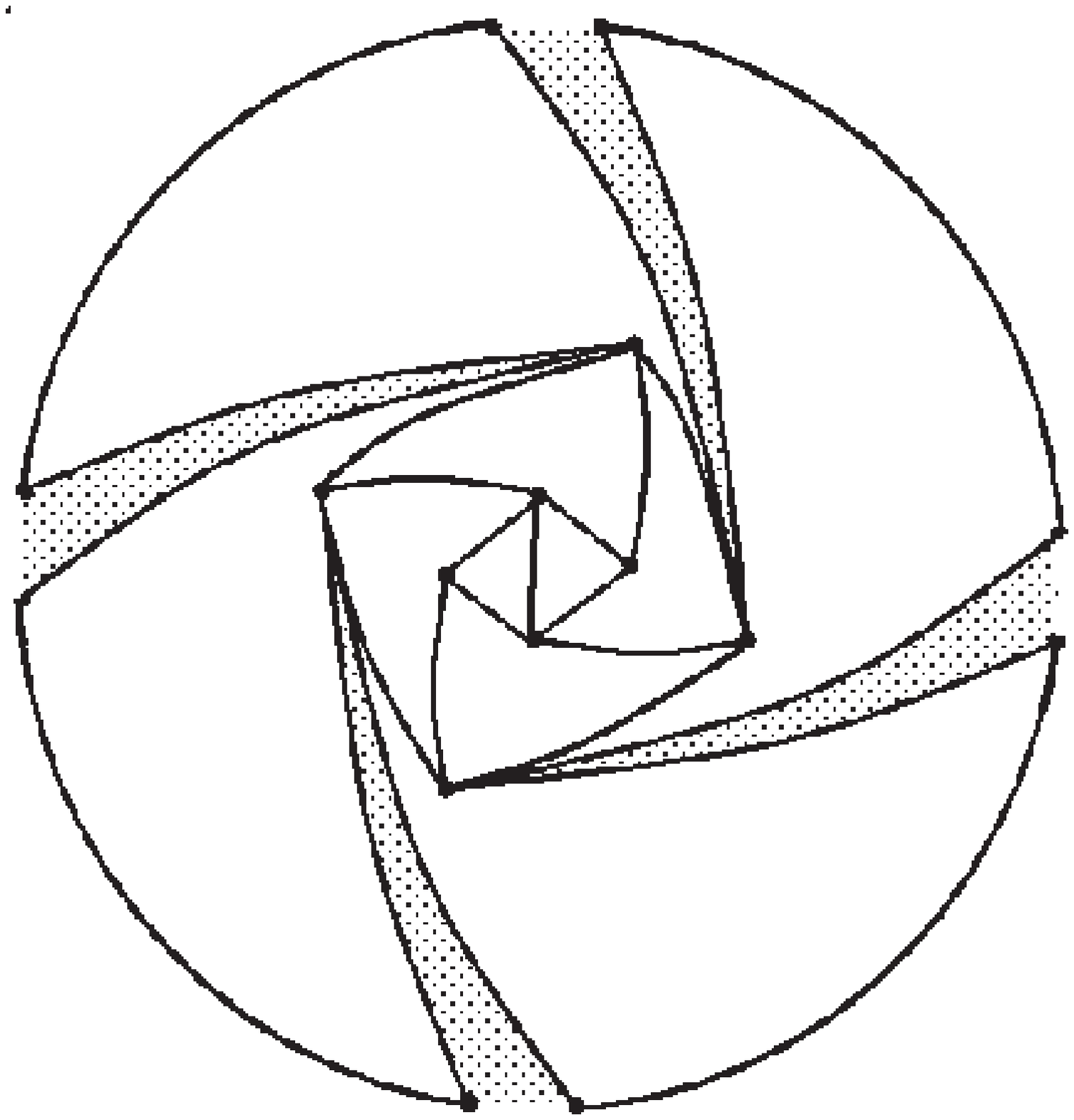}
\caption{}\label{}

\end{figure}

\textbf{A dictionary}

\bigskip

\begin{tabular}{|c|c|}
\hline
  Maxwell's reciprocals   of a spherical self-stressed graph  & virtual polytope
 \\\hline
Maxwell's reciprocals   of a \textbf{pointed} self-stressed graph
& \textbf{hyperbolic} virtual polytope \\

  \hline

  pointed pseudo-tiling of the 2-dimensional sphere $S^2$  & hyperbolic fan \\

  \hline
a pseudo-tiling of    $S^2$ with a non-zero self-stress  &
realizable hyperbolic fan
 \\\hline
  pseudo-di-gons of such a pseudo-tiling  & horns of the
hyperbolic polytope \\
\hline
3D lifting of a graph & graph of support function\\

\hline

negatively stressed edge of a  spherical graph &
blue edge of the hyperbolic fan \\

\hline

positively stressed edge of a  spherical graph &
red edge of the hyperbolic fan \\

\hline

\end{tabular}

\pagebreak

\textsc{REFERENCES}

\bigskip

[1]A.D.Alexandrov, \textit{On uniqueness theorems for closed
surfaces} (Russian), Doklady Akad. Nauk SSSR 22 No. 3(1939),
pp.99-102.


[2] R. Connelly, E. Demaine, G. Rote, \textit{Straightening
polygonal arcs and convexifying polygonal cycles}, Discrete Comput.
Geom., 30 (2003), pp. 205--239.

[3] H. Crapo, W. Whiteley,\textit{ Plane self-stresses and projected
polyhedra I: the basic pattern}, Structural Topology, 20 (1993), pp.
55--77.

[4] H. Crapo, W. Whiteley, \textit{Spaces of stresses, projections,
and parallel drawings for spherical polyhedra}, Beitrage zur Algebra
und Geometrie (Contributions to Algebra and Geometry), 35 (1994),pp.
259--281.

[5] R. Diestel, "Graph Theory", second ed., Springer-Verlag, Berlin,
2000.

[6] J. Graver, B. Servatius, H. Servatius, "Combinatorial Rigidity",
Graduate Studies in Mathematics,
 vol. 2, Amer. Math. Soc., 1993.

[7] R.Haas,  D. Orden,  G. Rote,  F Santos,
 B. Servatius,  H. Servatius, D. Souvaine,  I. Streinu, W. Whiteley,
 \textit{Planar minimally rigid graphs and pseudo-triangulations}, Comput.
Geom. 31 (2005), No.1-2, pp. 31-61.

[8] R. Haas, D. Orden, G. Rote, F. Santos, B. Servatius, H.
Servatius, D. Souvaine, I. Streinu, W. Whiteley, \textit{Planar
minimally rigid graphs}, in: Proc. 19th Ann. Symp. Comput. Geometry,
2003, pp. 154--163.

[9] T. Jordan, F. Zsolt,  W. Whiteley, \textit{An inductive
construction for plane laman graphs via vertex splitting}, S.
Albers, (ed.) et al., Algorithms - ESA 2004. 12th annual European
symposium, Bergen, Norway, September 14-17, 2004. Proceedings.
Berlin: Springer. Lecture Notes in Computer Science 3221(2004), pp.
299-310.

[10] A. Khovanskii, A. Pukhlikov, \textit{Finitely additive measures
of virtual polytopes},
  St. Petersburg Math. J. 4, No. 2 (1993),
  pp. 337-356.

[11] G. Laman, "On graphs and rigidity of plane skeletal
structures", J. Eng. Math. 4 (1970).

 [12] Y. Martinez-Maure, \textit{Contre-exemple \`a une caract\'erisation
conjectur\'ee de la sph\`ere}, C.R. Acad. Sci. Paris 332, No. 1
(2001), pp.41-44.

[13] Y.Martinez-Maure, \textit{Th\'eorie des h\'erissons et
polytopes}, C.R. Acad. Sci. Paris, Ser.1 336(2003), pp. 41-44.

[14]   P.McMullen, \textit{The polytope algebra},
 Adv.Math. 78, No. 1 (1989), pp. 76-130.

[15] D. Orden, G. Rote, F. Santos, B. Servatius, H. Servatius, W.
Whiteley, \textit{Non-crossing frameworks with non-crossing
reciprocals}, Discrete Comput. Geom. 32 No.4(2004), pp. 567-600 .

[16] D. Orden, F. Santos, \textit{The polytope of non-crossing
graphs on a planar point set}, Discrete Comput.
 Geom.,  http://arxiv.org/abs/math.CO/0302126.

[17] D. Orden, F. Santos, B. Servatius, H. Servatius,
\textit{Combinatorial pseudo-triangulations}, preprint, July 2003,
11 pages, http://arxiv.org/abs/math.CO/0307370.

[18] G.Panina,\textit{ New counterexamples to A.D. Alexandrov's
uniqueness hypothesis},
 Advances in Geometry, no. 5 (2005), pp. 301-317.

[19] G.Panina, \textit{A.D. Alexandrov's uniqueness theorem for
convex polytopes and its refinements }, to appear.

 [20] G.Panina, \textit{On hyperbolic virtual polytopes and hyperbolic
fans}, CESJM,  4 (2006), no 2, pp. 270-293.

[21] G.Panina, \textit{Isotopy problems for saddle  surfaces},
preprint of ESI,1796 (2006), available at
http://www.esi.ac.at/preprints/ESI-Preprints.html.

[22] M. Pocchiola, G. Vegter, \textit{Pseudo-triangulations: theory
and applications}, in: Proc. 12th Ann.
 Symp. Comput. Geometry, Philadelphia, 1996, pp. 291--300.

[23] M. Pocchiola, G. Vegter, \textit{Topologically sweeping
visibility complexes via pseudo-triangulations}, Discrete Comput.
Geom. 16 (1996), pp. 419-453.

[24] G. Rote, F. Santos, I. Streinu, \textit{Expansive motions and the
polytope of pointed pseudo
 triangulations}, in: B. Aronov, S. Basu, J. Pach, M. Sharir (Eds.), Discrete and Computational
 Geometry  - The Goodman-Pollack Festschrift, Springer-Verlag, Berlin, 2003, pp. 699--736.

[25] I. Streinu, \textit{A combinatorial approach for planar
non-colliding robot arm motion planning}, Proc. 41st Symp. Found. of
Comp.Sc.,Redondo Beach, California, 2000, pp. 443-453.

[26] I. Streinu,\textit{ Combinatorial roadmaps in configuration
spaces of simple planar polygons}, in: S. Basu, L. Gonzalez-Vega
(Eds.), Proc. DIMACS Workshop, Algorithmic and Quantitative Aspects
of Real Algebraic Geometry in Mathematics and Computer Science,
2003, pp. 181-206.

 [27] I. Streinu, \textit{Acute triangulations of polygons}, Discrete Comput. Geom. 34
(2005), no.4, 587-635.

[28] T.S. Tay, W. Whiteley, \textit{Generating isostatic graphs},
Structural Topology 11 (1985), pp. 21-68.


[29] W. Whiteley, \textit{Some matroids from discrete applied
geometry}, in: J. Bonin, J. Oxley,
 B. Servatius (Eds.), Matroid Theory, in: Contemporary Mathematics, vol. 197, Amer.
  Math. Soc., 1996, pp. 171-311.

\end{document}